# LOXODROMES AND GEODESICS ON ROTATIONAL SURFACES IN PSEUDO-ISOTROPIC SPACE

*Murat BABAARSLAN[1]*

[1]Yozgat Bozok University, Faculty of Science and Letters, Department of Mathematics, Yozgat, Türkiye.

[1]ORCID ID: https://orcid.org/0000-0002-2770-4126

[1]murat.babaarslan@bozok.edu.tr

*Ahmet SUNGUR[2]*

[2]Yozgat Bozok University, Graduate School of Natural and Applied Sciences, Department of Mathematics, Yozgat, Türkiye.

[2]ORCID ID: https://orcid.org/0000-0002-0569-7674

[2]ahmetsungur0001@gmail.com

**ABSTRACT**

Special curves and surfaces have an important place in mathematics, engineering and other fields of science. Loxodromes are special curves which cut all meridians on the Earth's surface at a constant angle and they are very popular in engineering. Ships sailing and airplanes flying along a fixed magnetic compass course move along this curve. The Mercator projections of the loxodromes on the sphere correspond to the lines and their stereographic projections to the logarithmic spirals. In general, loxodromes are not great circle arcs (geodesics). Geodesics correspond to the shortest distance between two points on the Earth's surface. Since loxodromes do not need to change course, they are important in navigation. Up till now, there have been many important studies of loxodromes on different surfaces (sphere, ellipsoid, rotational, helicoidal, canal, twisted, etc.) and in different ambient spaces (Euclidean, Minkowski, simply isotropic, etc.). However, there is no study on loxodromes in pseudo-isotropic space $\mathbb{I}_p^3$. In $\mathbb{I}_p^3$, different pseudo-isotropic angles can be defined depending on whether the vectors are space-like or time-like as in Minkowski space. In this study, we will first define these angles in $\mathbb{I}_p^3$. Then, the equations of space-like and time-like loxodromes and geodesics on rotational surfaces will be obtained. Some examples will also be provided.

**Keywords:** Loxodrome, geodesic, pseudo-isotropic space.

# 1. INTRODUCTION

The concept of a loxodrome on a rotational surface was given by Noble in 1905 [13]. Subsequently, there have been many studies of loxodromes on different surfaces (helicoidal, canal, twisted) and in different spaces (Euclidean, Minkowski) (see [1,4–9,12,15]). Loxodromes on invariant surfaces in three-manifolds were introduced by Caddeo et al. [3].

Loxodromes and geodesics on rotational surfaces in simply isotropic space were studied by Yoon in 2017 [16]. However, there is no study on loxodromes in $\mathbb{I}_p^3$. Therefore, in this study, we find space-like and time-like loxodromes and geodesics on rotational surfaces in $\mathbb{I}_p^3$ by using similar methods in [16].

# 2. PRELIMINARIES

In this section, we repeat some important notations in pseudo-isotropic space $\mathbb{I}_p^3$. For more details, we refer to [2,10,11].

The group of pseudo-isotropic motions in $\mathbb{I}_p^3$ is defined by

$$\begin{bmatrix} \bar{x} \\ \bar{y} \\ \bar{z} \end{bmatrix} = \begin{bmatrix} \cosh v & \sinh v & 0 \\ \sinh v & \cosh v & 0 \\ c_1 & c_2 & 1 \end{bmatrix} \begin{bmatrix} x \\ y \\ z \end{bmatrix} + \begin{bmatrix} a \\ b \\ c \end{bmatrix}, \tag{2.1}$$

where $a, b, c, c_1, c_2, v \in \mathbb{R}$.

On the $xy$−plane, $\mathbb{I}_p^3$ looks like Lorentzian plane $\mathbb{E}_1^2$. The projection of a vector $\boldsymbol{p} = (p_1, p_2, p_3) \in \mathbb{I}_p^3$ on $xy$−plane is called as the top view of $\boldsymbol{p}$ and it is given by $\widetilde{\boldsymbol{p}} = (p_1, p_2, 0)$. $z$−axis is invariant under the pseudo-isotropic motion given by (2.1) and a line with $z$−axis is called as an isotropic line and a plane which contains an isotropic line is called as an isotropic plane. A vector $\boldsymbol{p} = (0,0,p_3) \in \mathbb{I}_p^3$ is called as an isotropic vector and otherwise it is called as a non-isotropic vector.

Let $\boldsymbol{p} = (p_1, p_2, p_3)$ and $\boldsymbol{q} = (q_1, q_2, q_3)$ be vectors in $\mathbb{I}_p^3$. The pseudo-isotropic scalar product of $\boldsymbol{p}$ and $\boldsymbol{q}$ is given by

$$\langle \boldsymbol{p}, \boldsymbol{q} \rangle = \begin{cases} p_3 q_3 & \text{if } p_1 = p_2 = 0 \text{ and } q_1 = q_2 = 0 \\ p_1 q_1 - p_2 q_2 & \text{if } \quad \text{otherwise.} \end{cases} \tag{2.2}$$

A non-isotropic vector $\boldsymbol{p} \in \mathbb{I}_p^3$ is called as space-like, time-like and light-like if $\langle \boldsymbol{p}, \boldsymbol{p} \rangle > 0$ or $\boldsymbol{p} = \boldsymbol{0}$, $\langle \boldsymbol{p}, \boldsymbol{p} \rangle < 0$ and $\langle \boldsymbol{p}, \boldsymbol{p} \rangle = 0$ ($\boldsymbol{p} \neq \boldsymbol{0}$), respectively.

A curve $\alpha: I \to \mathbb{I}_p^3$ is called as regular if $\alpha'(t) \neq 0$ for all $t \in I \subset \mathbb{R}$. Also, it is called as admissible curve if all the osculating planes are not pseudo-isotropic. An admissible curve $\alpha: I \to \mathbb{I}_p^3$ is called as space-like, time-like or light-like if $\alpha'(t)$ is space-like, time-like or light-like, respectively.

Let $S: D \subset \mathbb{R}^2 \to \mathbb{I}_p^3$ be a smooth parametric surface in $\mathbb{I}_p^3$ with a coordinate system $\{u, v\}$, where $D$ is a open subset of $\mathbb{R}^2$. $S$ is called as admissible surface if its all tangent planes $T_p S$ are not pseudo-isotropic. It is easily seen that every admissible surface in $\mathbb{I}_p^3$ is time-like.

Now, we give the rotational surfaces in $\mathbb{I}_p^3$. Using the pseudo-isotropic motion which is given by (2.1), the hyperbolic rotation about $z$-axis in $\mathbb{I}_p^3$ is given by

$$\begin{bmatrix} \bar{x} \\ \bar{y} \\ \bar{z} \end{bmatrix} = \begin{bmatrix} \cosh v & \sinh v & 0 \\ \sinh v & \cosh v & 0 \\ 0 & 0 & 1 \end{bmatrix} \begin{bmatrix} x \\ y \\ z \end{bmatrix}, \qquad (2.3)$$

where $v \in \mathbb{R}$.

Let $\alpha_1(u) = (u, 0, f(u))$ be a space-like admissible curve lying in the isotropic $xz$-plane of $\mathbb{I}_p^3$, where $f$ is smooth function. If we rotate $\alpha_1$ (profile curve) about $z$-axis, we find the rotational surface with space-like meridian as follows:

$$R_1(u, v) = (u \cosh v, u \sinh v, f(u)). \qquad (2.4)$$

Similarly, let $\alpha_2(u) = (0, u, f(u))$ be a time-like admissible curve lying in the isotropic $yz$-plane of $\mathbb{I}_p^3$. If we rotate $\alpha_2$ (profile curve) about $z$-axis, we find the rotational surface with time-like meridian as follows:

$$R_2(u, v) = (u \sinh v, u \cosh v, f(u)). \qquad (2.5)$$

The first fundamental form (induced metric) of $R_1$ is

$$ds^2 = du^2 - u^2 dv^2 \qquad (2.6)$$

Similarly, the first fundamental form (induced metric) of $R_2$ is

$$ds^2 = -du^2 + u^2 dv^2 \qquad (2.7)$$

The pseudo-isotropic angles between two non-isotropic vectors in $\mathbb{I}_p^3$ can be defined as the Lorentzian angle between two vectors in $\mathbb{E}_1^2$ as follows.

Let $\boldsymbol{p}$ and $\boldsymbol{q}$ be space-like non-isotropic vectors in $\mathbb{I}_p^3$ that span a time-like vector subspace. Then, the pseudo-isotropic angle $\theta \in \mathbb{R}^+$ between $\boldsymbol{p}$ and $\boldsymbol{q}$ is given by

$$|\langle \boldsymbol{p}, \boldsymbol{q} \rangle| = \|\boldsymbol{p}\| \|\boldsymbol{q}\| \cosh \theta. \qquad (2.8)$$

Let $\boldsymbol{p}$ be space-like non-isotropic vector and $\boldsymbol{q}$ be time-like non-isotropic vector in $\mathbb{I}_p^3$. Then, the pseudo-isotropic angle $\eta \in \mathbb{R}^+ \cup \{0\}$ between $\boldsymbol{p}$ and $\boldsymbol{q}$ is given by

$$|\langle \boldsymbol{p}, \boldsymbol{q} \rangle| = \|\boldsymbol{p}\|\|\boldsymbol{q}\| \sinh \eta. \tag{2.9}$$

Let $\boldsymbol{p}$ and $\boldsymbol{q}$ be time-like non-isotropic vectors in $\mathbb{I}_p^3$ with same parity. Then, the pseudo-isotropic angle $\varphi \in \mathbb{R}^+ \cup \{0\}$ between $\boldsymbol{p}$ and $\boldsymbol{q}$ is given by

$$\langle \boldsymbol{p}, \boldsymbol{q} \rangle = \|\boldsymbol{p}\|\|\boldsymbol{q}\| \cosh \varphi \tag{2.10}$$

(see [14]).

A loxodrome in $\mathbb{I}_p^3$ is a curve which cuts all meridians of a rotational surface at a constant angle.

## 3. LOXODROMES ON ROTATIONAL SURFACES

In this section, we give the parametrizations of space-like and time-like loxodromes on rotational surfaces in $\mathbb{I}_p^3$.

We first can give the following theorem.

**Theorem 3.1.** The space-like loxodromes on the rotational surfaces with space-like meridians in $\mathbb{I}_p^3$ are

$$\gamma_1(t) = (\pm t\cosh\theta \cosh(\pm\tanh\theta \ln|t|), \pm t\cosh\theta \sinh(\pm\tanh\theta \ln|t|), f(\pm t\cosh\theta)), \tag{3.1}$$

where $\theta \in \mathbb{R}^+$.

**Proof.** We consider the meridian of the rotational surface $R_1$ given by (2.4) ($v = $ constant):

$$\beta_v(u) = (u\cosh v, u\sinh v, f(u)). \tag{3.2}$$

Then, the derivative of $\beta_v(u)$ is

$$\beta_v'(u) = (\cosh v, \sinh v, f'(u)). \tag{3.3}$$

We want to find a space-like loxodrome $\gamma_1$ on $R_1$ as follows:

$$\gamma_1(t) = (u(t)\cosh v(t), u(t)\sinh v(t), f(u(t))). \tag{3.4}$$

Assume that $\gamma_1$ is unit speed space-like curve, that is

$$\langle \dot\gamma_1(t), \dot\gamma_1(t) \rangle = \dot u^2(t) - u^2(t)\dot v^2(t) = +1. \tag{3.5}$$

If we take $\alpha_1(0) = \gamma_1(0) = (0,0,f(0))$, then we get $u(0) = 0$.

The pseudo-isotropic scalar product of $\dot{\gamma}_1(t)$ and $\beta'_{v}|_{\gamma_1(t)}$ is

$$\langle \dot{\gamma}_1(t), \beta'_{v}|_{\gamma_1(t)} \rangle = \dot{u}(t). \tag{3.6}$$

Let $\theta$ be a constant pseudo-isotropic angle given by (2.8) between $\gamma_1(t)$ and $\beta_{v}|_{\gamma_1(t)}$. Thus, we get

$$\dot{u}(t) = \pm\cosh\theta. \tag{3.7}$$

From (3.7), we have

$$u(t) = \pm t \cosh\theta. \tag{3.8}$$

Substituting (3.7) and (3.8) into (3.5), we have the following ordinary differantial equation (ODE):

$$\cosh^2\theta - t^2\cosh^2\theta \dot{v}^2(t) = 1. \tag{3.9}$$

The general solution of (3.9) is given by

$$v(t) = \pm\tanh\theta \ln|t|. \tag{3.10}$$

If we substitute (3.8) and (3.10) into (3.4), then we arrive to (3.1).

Thus, the proof is completed.

By using (2.9) and (2.10) and the similar techiques in the proof of Theorem 3.1, we can give the following theorems without proof.

**Theorem 3.2.** The time-like loxodromes on the rotational surfaces with space-like meridians in $\mathbb{I}_p^3$ are

$$\gamma_2(t) = (\pm t\sinh\eta \cosh(\pm\coth\eta \ln|t|), \pm t\sinh\eta \sinh(\pm\coth\eta \ln|t|), f(\pm t\sinh\eta)), \tag{3.11}$$

where $\eta \in \mathbb{R}^+ \cup \{0\}$.

**Theorem 3.3.** The space-like loxodromes on the rotational surfaces with time-like meridians in $\mathbb{I}_p^3$ are

$$\gamma_3(t) = (\pm t\sinh\eta \sinh(\pm\coth\eta \ln|t|), \pm t\sinh\eta \cosh(\pm\coth\eta \ln|t|), f(\pm t\sinh\eta)), \tag{3.12}$$

where $\eta \in \mathbb{R}^+ \cup \{0\}$.

**Theorem 3.4.** The time-like loxodromes on the rotational surfaces with time-like meridians in $\mathbb{I}_p^3$ are

$$\gamma_4(t) = (-t\cosh\varphi \sinh(\pm\tanh\varphi \ln|t|), -t\cosh\varphi \cosh(\pm\tanh\eta \ln|t|), f(-t\cosh\varphi)), \quad (3.13)$$

where $\varphi \in \mathbb{R}^+ \cup \{0\}$.

Finally, we give the following example.

**Example 3.5.** We consider a rotational surface with space-like meridian obtained by $\alpha_1(u) = (u, 0, e^u)$. If we take $\theta = \pi/4$ and $t \in (1,2)$, we can plot the space-like loxodrome (green) and the meridian (blue) in Figure 1 as follows.

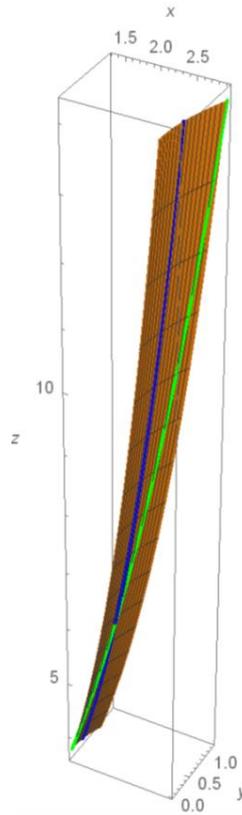

**Figure 1.** Rotational surface with space-like meridian (orange); space-like loxodrome (green), meridian (blue).

## 4. GEODESICS ON ROTATIONAL SURFACES

In this section, we give the parametrizations of geodesics on rotational surfaces in $\mathbb{I}_p^3$.

**Theorem 4.1.** Let $\gamma_5(t) = (u(t)\cosh v(t), u(t)\sinh v(t), f(u(t)))$ be a curve on the rotational surface with space-like meridian and $\gamma_6(t) = (u(t)\sinh v(t), u(t)\cosh v(t), f(u(t)))$ be a curve on the on the rotational surface with time-like meridian in $\mathbb{I}_p^3$. Then,

**(i)** The parallels of the rotational surfaces ($u =$ constant) are geodesics iff $\dot{v}(t) = 0$, when $t = t_0$.

**(ii)** The meridians of the rotational surfaces ($v =$ constant) are geodesics iff $u(t) = at + b$ for $a \neq 0, \ b \in \mathbb{R}$.

**(iii)** The curves $\gamma_5(t)$ and $\gamma_6(t)$ which are given by

$$u(t) = \mp \frac{1}{\sqrt{|c_1|}}\sqrt{(c_1 t + c_2)^2 - c^2}, \quad v(t) = \frac{1}{2}\ln\left|\frac{c_1 t + c_3}{c_1 t + c_4}\right| + c_5$$

are geodesics, where $c, c_1 \neq 0, c_2, \ c_3 = (c_2 - c), \ c_4 = (c_2 + c), c_5 \in \mathbb{R}$.

**Proof.** Let $\gamma_5(t) = R_1(u(t), v(t))$ be a geodesic of the rotational surface with space-like meridian in $\mathbb{I}_p^3$. The Euler Lagrange equations for (2.6) are

$$\begin{cases} \frac{d}{dt}(\dot{u}(t)) = -u(t)\dot{v}^2(t), \\ \frac{d}{dt}(u^2(t)\dot{v}(t)) = 0. \end{cases} \quad (4.1)$$

We consider a parallel of the rotational surface with space-like meridian ($u =$ constant). From the first equation of (4.1), we have $\dot{v}(t) = 0$. Thus, the second equation of (4.1) is also satisfied.

We consider a meridian of the rotational surface with space-like meridian ($v =$ constant). From the first equation of (4.1), we have $u(t) = at + b$ for $a \neq 0, \ b \in \mathbb{R}$.

We now investigate any geodesics of $R_1$. From the second equation of (4.1), we get

$$\dot{v}(t) = \frac{c}{u^2(t)}, \quad (4.2)$$

where $c$ is a non-zero constant. If we substitute (4.2) into the first equation of (4.1), we get

$$u^3(t)\ddot{u}(t) = -c^2. \quad (4.3)$$

If (4.3) is multiplied by $2\dot{u}(t)\frac{1}{u^3(t)}$ and the obtained equation is integrated, we find

$$\dot{u}^2(t) = \frac{c_1 u^2(t) + c^2}{u^2(t)}, \tag{4.4}$$

where $c_1$ is a non-zero constant.

The general solution of (4.4) is

$$u(t) = \mp \frac{1}{\sqrt{|c_1|}}\sqrt{(c_1 t + c_2)^2 - c^2}, \tag{4.5}$$

where $c_2$ is a constant. Substituting (4.5) into (4.2), we get

$$v(t) = \frac{1}{2}\ln\left|\frac{c_1 t + c_3}{c_1 t + c_4}\right| + c_5, \tag{4.6}$$

where $c_3 = (c_2 - c)$, $c_4 = (c_2 + c)$ and $c_5$ are constants.

If we make the similar computations for the curve $\gamma_6(t) = (u(t)\sinh v(t), u(t)\cosh v(t), f(u(t)))$ on the rotational surface with time-like meridian given by (2.5), we arrive to the same equations.

Thus, the poof is completed.

Finally, we give the following example.

**Example 4.2.** We consider a rotational surface with time-like meridian obtained by $\alpha_2(u) = (0, u, \cos u)$. If we take $c = 1, c_1 = 4, c_2 = 2, c_5 = 0, a = 2, b = 5$ and $t \in (0,2)$, we can plot the three cases from Theorem 4.1 in Figure 2 as follows.

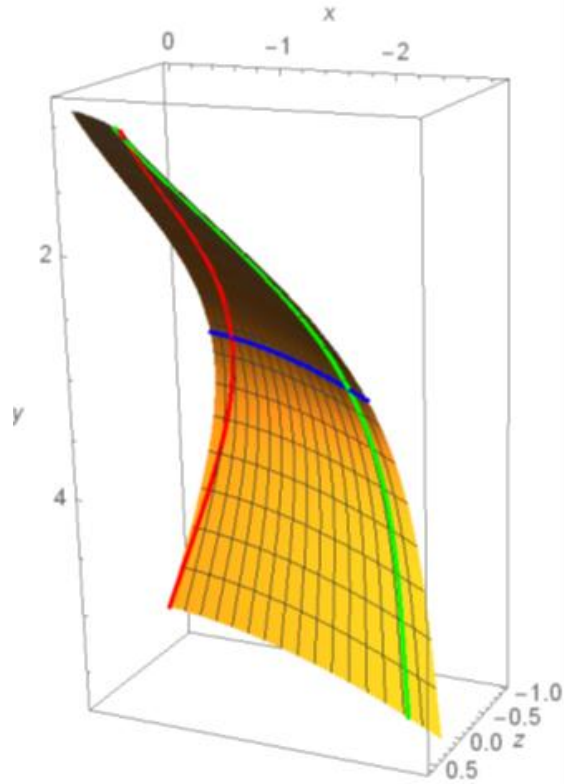

**Figure 2.** Rotational surface with time-like meridian (orange); parallel (blue), meridian (green), another meridian (red).

## ACKNOWLEDGEMENT

This work is a part of the master thesis of the second author.